\documentclass[reqno]{amsart}
\usepackage{amsfonts}
\usepackage{amsmath}
\usepackage{amsthm}
\usepackage{amssymb}
\usepackage{oldgerm}
\usepackage{fancyhdr}
\usepackage{fancybox}
\usepackage{mathrsfs}
\usepackage{epsfig}
\def\phi{\varphi }
\swapnumbers

\theoremstyle{plain}
\newtheorem{theorem}{Theorem}[section]

\newtheorem{lemma}[theorem]{Lemma}
\newtheorem{proposition}[theorem]{Proposition}
\theoremstyle{definition}
\newtheorem{definition}[theorem]{Definition}
\theoremstyle{remark}

\newtheorem*{remark}{Remark}

\newtheorem{example}[theorem]{Example}
\numberwithin{equation}{section}

\newcommand{\lieg}{\mathfrak{g}}

\newcommand{\lieu}{\mathfrak{u}}
\newcommand{\liek}{\mathfrak{k}}
\newcommand{\liep}{\mathfrak{p}}
\newcommand{\lieq}{\mathfrak{q}}
\newcommand{\liea}{\mathfrak{a}}
\newcommand{\lieb}{\mathfrak{b}}

\newcommand{\xdiag}{\underline x}
\newcommand{\FF}{\mathbb{F}}
\newcommand{\CC}{\mathbb{C}}
\newcommand{\RR}{\mathbb{R}}
\newcommand{\HH}{\mathbb{H}}
\newcommand{\ZZ}{\mathbb{Z}}
\newcommand{\NN}{\mathbb{N}}
\newcommand{\singspec}{\textnormal{spec}_s}

\newcommand{\Waff}{W_{\textnormal{aff}}}

\newcommand{\la}{\lambda_{\alpha}}
\renewcommand{\Re}{\textnormal{Re} \,}

\title{Convolution algebras for Heckman-Opdam polynomials derived from compact Grassmannians}
\author{Heiko Remling}
\address{H. Remling:}
\email{hremling@gmail.com}
\author{Margit R\"osler}                   
\address{M. R\"osler (corresponding author): Institut f\"ur Mathematik, Universit\"at Paderborn,
Warburger Str. 100, D-33098 Paderborn, Germany}
\email{roesler@math.upb.de}


\begin{document}

\date{}

\begin{abstract}
We study convolution algebras associated with Heckman-Opdam polynomials. For root systems of type $BC$ we 
derive three continuous classes of positive convolution algebras (hypergroups) by interpolating the double coset convolution structures of compact 
Grassmannians $U/K$ with fixed rank over the real, complex or quaternionic numbers. These convolution algebras are linked to explicit positive product 
formulas for  Heckman-Opdam polynomials of type $BC$,
 which occur  for certain discrete multiplicities as the spherical functions of $U/K$. The results complement those of \cite{Ro09}
for the non-compact case. 
\end{abstract}

\maketitle

%
\section{Introduction}
%

In the theory of multivariable hypergeometric functions and polynomials of Heckman, Cherednik and Opdam, 
the existence of product formulas and positive convolution algebras is in general unsolved.
In \cite{Ro09}, three continuous series of positive convolution algebras having Heckman-Opdam hypergeometric functions  as multiplicative 
functions were obtained by interpolating geometric cases in an explicit way,  namely the product formulas for the spherical functions of non-compact 
Grassmannians. In these cases, a full picture of harmonic analysis for the hypergeometric transform could be obtained. The present paper extends these 
results to the dual situation related to compact Grassmannians and convolution algebras for three continuous series of Heckman-Opdam Jacobi polynomials of 
type $BC$. The Heckman-Opdam Jacobi polynomials occur as spherical functions of compact Grassmannians. This observation was first made by Koornwinder in 
rank two (see \cite{Koorn2}), and the corresponding two-variable analogues of Jacobi polynomials were  introduced by Koornwinder in \cite{Koorn1}.
There is a broad literature on multivariable Jacobi polynomials, in particular in the context of multivariate statistics. We mention \cite{JC} and \cite{La}. However, 
explicit product formulas 
have not been given so far.

In the present paper we start from the compact Grassmannians $G_{p+q,q}(\mathbb
F)$ of
$p$-dimensional subspaces of $\mathbb F^{p+q}$ where $p>q$ and $\mathbb F$ is
one of the (skew) fields $\mathbb F= \mathbb R, \mathbb C, \mathbb H$. The
Grassmanians $G_{p+q,q}(\mathbb F)$ are realized as a Riemannian symmetric
spaces $U/K$ with $U=SU(p+q,\mathbb F)$ and $K=S(U(p,\mathbb F)\times
U(q,\mathbb F))$. 
They are  dual to the non-compact Grassmannians studied in \cite{Ro09}. Similar as in loc.cit., we write down the product formula for their spherical functions in an explicit way which allows analytic continuation with respect to the dimension parameter $p$, the rank $q$ being fixed. 
The spherical functions in the geometric cases (corresponding to integral $p$) are Heckman-Opdam Jacobi polynomials of type $BC_q$ 
(for $\mathbb F= \mathbb C, \,\mathbb H$) or $B_q$ (for $\mathbb F=\mathbb R),$
with certain discrete multiplicities. In contrast to the non-compact case, the
study of the geometric background needs some care in the compact case,
especially for $\mathbb F= \mathbb R.$

Our continuation gives an explicit  product formula for an interpolated
continuous range of multiplicities.
This formula  in part generalizes Koornwinder's product formula for Jacobi
polynomials \cite{Ko74}  to higher rank. Naturally, it is  similar to the
non-compact case, but  direct analytic continuation from the non-compact to the
compact case seems  not feasible.
We obtain three continuous classes of commutative hypergroup algebras on the fundamental alcove of the associated
affine reflection group, with the associated Heckman-Opdam Jacobi polynomials as characters. 

The organisation of this paper is as follows:
In Section \ref{sec_dunkl} we recall some basics of trigonometric Dunkl theory.
Section \ref{Symm_Raum} is a summary of the necessary background from the
theory of compact symmetric spaces. After that, we start in Section 4 with the
compact
Grassmannians  $U/K$, identify their spherical functions with
Heckman-Opdam Jacobi
polynomials, and use a Cartan decomposition of $U$  to make their
product formula
explicit.  Following the idea of \cite{Ro09},  this product formula is then
analytically continued. Section 5 contains a review of the rank one case. In
Section 6, the related hypergroup structures and their dual spaces on the
fundamental alcove are studied.

%
\section{Fundamentals of Trigonometric Dunkl Theory}\label{sec_dunkl}
%

This section is a 
short review of the fundamentals of trigonometric Dunkl theory which will be needed in this article. For details, we refer to the work of Heckman and Opdam 
(\cite{HS94}, \cite{Op95}, \cite{Op00}). \\

Let $\liea$ be a $q$-dimensional Euclidean space with inner product $\langle
\cdot, \cdot \rangle$
which is extended to a complex bilinear form on the 
complexification $\liea_{\CC}$ of $\liea$. 
We identify $\liea$ with its dual space $\liea^*$ via the given inner product.
Let $\Sigma \subset \liea$ be a
(not necessarily reduced) root system.
For $\alpha \in \Sigma$ we write $\alpha^{\vee} := 2 \alpha/\langle \alpha, \alpha \rangle$ for the coroot of $\alpha$ and denote by 
$s_{\alpha} (x) = x - \langle \alpha^{\vee}, x \rangle \alpha$ the reflection in the hyperplane $H_{\alpha}$ perpendicular to $\alpha$.

The reflections $\{ s_{\alpha} \, : \, \alpha \in \Sigma\}$ generate the Weyl
group $W=W(\Sigma)$.
We define the root lattice $Q := \ZZ.\Sigma$ 
and the coroot lattice $Q^{\vee} = \ZZ.\Sigma^{\vee}$. We fix some
positive subsystem  $\Sigma^+$ of $\Sigma$, and write $\liea^+ := \{ \lambda \in
\liea \, : \, \langle \lambda, \alpha \rangle > 0 \, \,  \forall \alpha \in
\Sigma^+ \}$ for the associated Weyl chamber.
For $\alpha \in \Sigma$ and $\lambda \in \liea_{\CC}$ let
\begin{equation*}\label{lambda_alpha}
\lambda_{\alpha} := \frac{ \langle \lambda, \alpha \rangle}{\langle \alpha, \alpha \rangle}.
\end{equation*}
The weight lattice associated with $\Sigma$ is given by
\[
P = P(\Sigma) := \{ \lambda \in \liea \, : \, \lambda_{\alpha} \in \ZZ \,\textnormal{ for all } \alpha \in \Sigma \}
\]
and the cone of dominant weights associated with $\Sigma^+$ is
\[
P^+ = P^+(\Sigma):= \{ \lambda \in \liea \, : \, \lambda_{\alpha} \in \ZZ^+ \,\textnormal{ for all } \alpha \in \Sigma^+ \}
\]
Here $\ZZ^+ := \{0,1,2, \ldots \}$. Notice that $2\Sigma\subset P.$ The cone  $Q^+ = \ZZ^+. \Sigma^+$ defines a natural partial ordering $\preceq$ on $P^+$:
\begin{equation*}
\mu \preceq \lambda \iff \lambda - \mu \in 2Q^+.
\end{equation*}





A multiplicity function on $\Sigma$ is a $W$-invariant map
$m: \Sigma \to \CC$, $\alpha \mapsto m_{\alpha}$. We denote the space of
multiplicity functions on 
$\Sigma$ by $\mathcal M$ and define
\begin{equation}\label{rho}
\rho = \rho(m) := \frac{1}{2} \sum_{\alpha \in \Sigma^+} m_{\alpha} \alpha.
\end{equation}
\begin{definition}\label{Def_DO}
Let $\xi \in \liea$ and $m\in \mathcal M.$ The Dunkl-Cherednik operator associated with $\Sigma$ and $m$ is given by
\[
T_{\xi} = T(\xi,m) := \partial_{\xi} + \sum_{\alpha \in \Sigma^+} m_{\alpha} \langle \alpha, \xi \rangle \frac{1}{1-e^{-2 \alpha}} (1- s_{\alpha}) - \langle \rho , \xi \rangle,
\]
where $\partial_{\xi}$ is the usual directional derivative and $e^{\lambda}(z)
:= e^{\langle \lambda, z \rangle}$ for $\lambda, z \in \liea_{\CC}$.
\end{definition}

\begin{remark}\label{Bem_Not}
Heckman and Opdam use a slightly different notation. They
consider a root system $R$ with  multiplicity $k$, which is connected to our
notation via
\[
R = 2 \Sigma, \quad k_{2 \alpha} = \frac{1}{2} m_{\alpha}.
\]
Note that this implies further differences. Our notation comes from the theory of symmetric spaces. 
\end{remark}

For fixed multiplicity $m$, the operators $T_{\xi}$, $\xi \in \liea_{\CC}$ commute.
Therefore the assignment $\xi \mapsto T(\xi,m)$ uniquely extends to a homomorphism on the symmetric algebra $S(\liea_{\CC})$ over $\liea_{\CC}$, which may be identified with the algebra of complex polynomials on $\liea_{\CC}$. Let $T(p,m)$ be the operator which corresponds in this way to $p \in S(\liea_{\CC})$. If $p\in   S(\liea_{\CC})^W$, the subspace of $W$-invariant polynomials  on $\liea_{\CC}$, 
 then $T(p,m)$ acts as a differential operator on the space of $W$-invariant analytic functions on $ \liea$.
Consider the so-called hypergeometric system
\[ T(p,m) \phi = p(\lambda) \phi  \quad \text{ for all }\, p \in S(\liea_{\CC})^W\]
where $\lambda\in \liea_{\CC}$ is a fixed spectral parameter. According to fundamental results by Heckman and Opdam (see \cite{HS94}), 
there exists an open set of regular multiplicities $\mathcal M^{reg}\subset \mathcal M$, containing all nonnegative multiplicities $m\geq 0$, 
such that for each fixed spectral parameter $\lambda$ and each $m\in \mathcal M^{reg}$, the associated hpergeometric system
has a unique $W$-invariant solution $\phi =F_{\lambda}(m;\cdot) = F(\lambda,m; \cdot)$ which is analytic on $\liea$ and satisfies $F_\lambda(m;0)=1$.
Moreover, there is a $W$-invariant tubular neighborhood $U$ of $\liea$ in $\liea_{\CC}$ such that $F$ extends to a (single-valued) holomorphic function 
$F: \liea_{\CC} \times  \mathcal M^{\mathrm{reg}}\times U \to \CC$. The function $F(\lambda,m;x)$ 
is $W$-invariant in both $\lambda$ and $x$. It is called the hypergeometric function associated with $\Sigma$. 
For certain spectral parameters $\lambda$, the functions $F_\lambda$ are
actually exponential polynomials, the so-called Heckman-Opdam polynomials.
In order to make this precise, we need some more notation. 

Let $\mathcal T := \text{lin} \{e^{i\lambda} \, : \, \lambda \in P \}$ be the
space of trigonometric polynomials associated with $P$. Trigonometric
polynomials are $\pi  Q^{\vee}$-periodic. Let
\[ M_\lambda := \sum_{\mu \in W.\lambda} e^{i\mu}, \quad \lambda \in P^+ \]
denote the $W$-invariant orbit sums. They
form  a basis of the space of
$W$-invariant trigonometric polynomials $\mathcal T^W$. 

For nonnegative multiplicity $m$, consider the
$W$-invariant weight function
\begin{equation}
w_m(x):= \prod_{\alpha \in \Sigma^+} \left| e^{i\langle\alpha,x\rangle} - e^{- i\langle \alpha,x\rangle} \right|^{m_{\alpha}}.
\end{equation}
on the torus $T:= \liea/\pi Q^\vee.$ The Heckman-Opdam polynomials associated
with $\Sigma$ and $m\geq 0$ are defined by
\[
 P_{\lambda} =P_\lambda(m; \cdot) := \sum_{\mu \in P^+, \, \mu \preceq \lambda}
c_{\lambda \mu}(m) M_\mu , \quad \lambda \in P^+
\]
where the coefficients $c_{\lambda \mu}(m)$ are uniquely determined by the conditions 
\begin{enumerate}\itemsep=-1pt
 \item[\rm{(i)}] $c_{\lambda \lambda}(m) = 1$
 \item[\rm{(ii)}] $P_\lambda$ is orthogonal to $M_\mu$ in $L^2(T; w_m)$  for all $\mu \in P^+$ with $ \mu \prec \lambda$. 
\end{enumerate}

\begin{remark}\label{Bem_i}
Notice that our notion slightly differs from that of Heckman and Opdam (e.g.
\cite{HS94}, \cite{Op00}),
namely by a factor $i$ in the spectral variable. This choice of notation will be
more convenient for our purposes.
\end{remark}

The polynomials $P_\lambda$  form an orthogonal basis of $L^2(T,w_m)^W$,
the subspace of $W$-invariant elements from $L^2(T,w_m)$. Their coefficients
$c_{\lambda \mu}(m)$ are rational functions in the
$m_\alpha.$ Moreover, their numerator
and denominator polynomials have nonnegative integral coefficients. This was
observed in \cite{Mac87}, Par.11. As a consequence,
\begin{equation}\label{conjugation} P_\lambda(-z) =
\overline{P_\lambda(\overline z)} \quad\text{for all }z\in
\liea_\mathbb C,\end{equation}
c.f. \cite{RR}. Moreover, 
the function $(m, z) \mapsto
P_\lambda(m;z)$
uniquely extends to a holomorphic function on $\{m\in \mathcal M: \text{Re}\, m
>0 \}\times \liea_{\mathbb C}.$

The connection between the Heckman-Opdam polynomials and the hypergeometric function is as follows (see \cite{HS94}, Section 4.4):

\begin{lemma}\label{Hyper_Jacobi} 
For all  $\lambda \in P^+$ and $m\geq 0$, the function
$F_{\lambda+\rho}(m ;\,\cdot\,)$ extends holomorphically to $\liea_\mathbb C$
with 
\[
F_{\lambda+\rho}(m ; iz) = c(\lambda + \rho,m) P_\lambda(m ; z).
\]
Here  the $c$-function $c(\lambda + \rho,m)$ is given by
\[
c(\lambda + \rho,m) = \prod_{\alpha \in \Sigma^+} \frac{\Gamma(\la + \rho_{\alpha} + \frac{1}{4} m_{\alpha/2}) \Gamma(\rho_{\alpha} + \frac{1}{4} m_{\alpha/2} + \frac{1}{2} m_{\alpha})}{ \Gamma(\lambda_{\alpha} + \rho_{\alpha} + \frac{1}{4} m_{\alpha/2} + \frac{1}{2} m_{\alpha}) \Gamma(\rho_{\alpha} + \frac{1}{4} m_{\alpha/2})}.
\]
with the convention that $m_{\alpha/2} = 0$ if $\alpha/2\notin \Sigma.$ 
\end{lemma}

We shall work with a renormalized version of the Heckman-Opdam polynomials, defined by
\begin{equation}\label{Def_R}
 R_{\lambda} (z) :=R_\lambda(m;z) :=  c(\lambda + \rho,m) P_{\lambda}(m;z) =
F_{\lambda+\rho}(m ; iz).
\end{equation}
Thus $R_{\lambda} (0) = 1.$  

The periodicity and $W$-invariance of the Heckman-Opdam polynomials is described by the affine Weyl group
\[ \Waff = \pi Q^{\vee} \rtimes W.\]
This is the Coxeter group generated by the affine reflections in the hyperplanes
\[ H_{\alpha, k} := \{ x  \in \liea : \langle \alpha, x \rangle = k \}= \,H_{\alpha,0} + \frac{k}{2}\alpha^\vee, 
\,\, k\in \pi\mathbb Z,\, \alpha \in \Sigma.\]
A fundamental domain for the action of $W_{\text{aff}}$ on $\liea$ is given by the closed fundamental alcove
\[A_0 = \{ x \in \liea : 0 \le \langle \alpha, x \rangle \le \pi \,\,\textnormal{ for all } \,\alpha \in \Sigma^+\}.\]
The trigonometric polynomials from $\mathcal T^W$ are $\Waff$-invariant, and can therefore be considered
as functions on the alcove $A_0$. Note that the Jacobi polynomials $R_{\lambda}\,, \lambda\in P^+$ form an orthogonal basis of $L^2(A_0, w_m).$

%
\section{Compact symmetric spaces and their spherical functions}\label{Symm_Raum}
%

In this section we recall some general background from the theory of symmetric
spaces. Standard references are the monographs \cite{Hel78}, \cite{Hel84},
\cite{T}.

Let $U/K$ be a Riemannian symmetric space of the compact type, where
$U$ is a connected compact Lie group and $K$ is a closed subgroup
such there exists an involutive automorphism 
$\theta: U \to U$ with $U^{\theta}_0  \subseteq K \subseteq U^{\theta}$. Here
$U^{\theta} = \{ u \in U \, : \, \theta (u) = u \}$ and $U^{\theta}_0$
denotes the identity component of $U^{\theta}$. To avoid technicalities, we
assume in the following that  $U$ is semisimple and $K$ is connected. Note that
if $U$ is simply connected, then $U^{\theta}_0$ is connected and $K=
U^{\theta}_0.$

The derivation of $\theta$ gives an involution of the Lie algebra $\lieu$ of
$U$. We write the associated Cartan decomposition of $\lieu$ as
  $\lieu = \liek \oplus \lieq$ with $\liek = \{ X\in \lieu: \theta(X) = X\}, \,
\lieq = \{ X\in \lieu: \theta(X) =-X\}.$ Let $\lieb\subset \lieq$ be a maximal
abelian subspace, and put $\liep = i\lieq.$ Then $\liea := i\lieb$ is a maximal
abelian subspace of $\liep$. 
Denote by $G$ the connected real Lie subgroup of the complexification $U_{\CC}$
of $U$ with Lie algebra
$\lieg := \liek \oplus \liep$. $G$ is a noncompact semisimple Lie
group with $K\subseteq G$, and $G/K$ is a symmetric space of the
non-compact type, called the
non-compact dual of $U/K.$ A Cartan involution $\tau$ of $G$ with $K= G^\tau$
is given by $\tau = \theta_{\mathbb C}\vert_G$, where $\theta_\mathbb C$ is the
analytic 
continuation of $\theta$ to $U_\mathbb C$.  Let $\liea$ be a maximal abelian
subspace of $\liep$.  Then $\liea$ is a finite-dimensional Euclidean space with
the Killing form
$B(\,.\,,\,.\,)$ as scalar product.
We shall identify $\liea$ with its dual $\liea^*$ via
$B$.
Further, we denote by $\Sigma := \Sigma (\lieg, \liea)$ the restricted root
system of
$\lieg$ with respect to $\liea$ 
 and by $\Sigma^+$ a fixed subset of positive
restricted roots. 

Recall that for an arbitrary Lie group $G$ with compact subgroup $K$, a spherical function of $(G,K)$ is a 
nonzero, $K$-biinvariant function $\varphi: G \to \CC$ which satisfies the product formula
\begin{equation}\label{prodformula}
\phi(g)\phi(h) = \int_K \phi(gkh)dk \quad \text{ for all }\, g,h\in G
\end{equation}
where $dk$ denotes the normalized Haar measure on $K$.

Assume now that $(U,K)$ and $(G,K)$ are as above. The spherical functions of
$(G,K)$ are given by the
the Harish-Chandra formula
\[ \varphi_\lambda(g) = \int_K e^{B(\lambda-\rho, H(gk))} dk, \quad
\lambda \in \liea_{\mathbb C}.\]
Here $\rho = \frac{1}{2}\sum_{\alpha\in \Sigma_+} m_\alpha \alpha$ where
 $m_\alpha$ is the multiplicity (that is,  the dimension) of the root space
associated with $\alpha$, and for $g\in G, \,H(g)\in A$ denotes the unique
abelian part of $g$ in the
Iwasawa decomposition $G=KAN$ associated with $\Sigma^+$. We have
$\varphi_\lambda = \varphi_\mu$ iff the orbits of $\lambda$ and $\mu$ under the
Weyl group $W(\Sigma)$ coincide.

The spherical functions of $(U,K)$ are all positive-definite and are obtained as
matrix coefficients of the $K$-spherical irreducible representations of $U$
(\cite{Hel84}, Chap. IV,  Theorems 3.4 and 4.2).
Recall that an irreducible unitary representation $\pi$ of $U$ in a Hilbert
space $V$ is called $K$-spherical if
the space \[V^K = \{ v\in V: \pi(k)v = v \,\text{ for all }\, k\in K\}\] of
$K$-fixed
vectors is different from $\{0\};$ in this case, actually $\dim V^K = 1,$
because $(U,K)$ is a Gelfand pair.
The spherical functions of $(U,K)$ are parametrized by the set $\Lambda_K(U)$
of (restrictions of) highest weights of $K$-spherical irreducible
representations of $U.$ If $U$ is simply connected, then by the Cartan-Helgason
theorem (\cite{Hel84}, Theorem 4.1, Chap. V), $\Lambda_K(U)$ coincides with the
set $P^+(\Sigma)$. This remains true if not $U$, but $U/K$ is simply connected,
see \cite{T}, Theorem 8.2. and Corollary 1, or \cite{OP}, Section 1.3.
The spherical function associated with $\mu\in \Lambda_K(U)$ is given by
\[ \psi_\mu(u) = \langle \pi_\mu(u)e_\mu, e_\mu\rangle\]
where $(\pi_\mu, V_\mu)$ is the spherical representation associated with $\mu$
and $e_\mu\in V_\mu^K$  is  a $K$-fixed vector with $\|e_\mu\|=1.$

We assume again that $U$ or $U/K$ is simply connected, $K$ is
connected, and $G/K$ is
the non-compact dual of $U/K$. Then there is the following  close
connection between the
spherical functions on $G$ and those on $U$.

\begin{proposition}\label{Sph_Dual}
Every spherical function $\phi_{\mu}$ of $(G,K)$ $(\mu \in \liea_{\CC})$ is 
analytic on $G$. It extends to a holomorphic function on the
complexification $G_{\mathbb C}= U_{\mathbb C}$ if and only if $\mu$ is
contained in the $W(\Sigma)$-orbit of $\lambda + \rho$ for some $\lambda \in
P^+(\Sigma)$. In this case, we denote the analytic extension also by
$\phi_{\mu}$. The restriction of this extension to $U$ is a spherical
function on $U.$ More precisely, we have the identity
\[
\phi_{\lambda + \rho} \vert_U = \psi_{\lambda}, \qquad \lambda \in P^+(\Sigma).
\]
Conversely, each spherical function $\psi_{\lambda}$ of $(U,K)$
extends to a holomorphic function  $\psi_{\lambda}$ on
$U_{\CC}$ and its restriction to $G$ coincides with the spherical function
$\varphi_{\lambda + \rho} $ of $(G,K)$.
\end{proposition}

\begin{proof}
\cite{Hel84}, Chap.V, proof of Theorem 4.4, and Lemma 2.5 in \cite{BOP05}.
\end{proof}

The following important fact links the theory of Heckman and
Opdam with the classical theory of symmetric spaces.

\begin{proposition}\label{spherical_hyper} (\cite{HS94}, Theorem 5.2.2)
Let $\phi_{\mu}$, $\mu \in \liea_{\CC}$, be a spherical function of $(G,K)$,
and let $\Sigma$ and $m$ be the associated restricted root system and geometric
multiplicity.
Then for $x\in \liea, $
\[
\phi_{\mu}(\exp x)   = F_\mu(m;x).
\]
\end{proposition}

Combining this with Proposition \ref{Sph_Dual} and
with  \eqref{Def_R}, we obtain

\begin{theorem}\label{ident_C_H}
The spherical
functions of $(U,K)$  - restricted to $\exp (i\liea)$ - are Heckman-Opdam
polynomials of type $\Sigma$ and with multiplicity $m$: For all $x\in \liea$ and
$ \lambda \in P^+(\Sigma),$
\[
\psi_{\lambda} (\exp (ix)) = \phi_{\lambda + \rho} (\exp (ix)) =
F_{\lambda+\rho} (m;ix) =  R_{\lambda}(m;x).
\]
\end{theorem}

The second equality in the theorem above follows from Proposition
\ref{spherical_hyper}
since $\phi_{\lambda + \rho}$ and $F_{\lambda + \rho}$ are holomorphic on
$G_\mathbb C$ and $\liea_\mathbb C$, respectively.

%
\section{A product formula for Heckman-Opdam polynomials of type $BC$}
%

Let $\mathbb F$ be one of the (skew) fields $\mathbb R, \mathbb C, \mathbb H$
with
the standard involution
 $x\mapsto \overline x$ and norm $|x| = (\overline x x)^{1/2}$. By
$M_{n}(\mathbb F)$ we denote the set of $n\times n$ matrices over $\mathbb
F$, also viewed as $\mathbb F$-linear transformations from $\mathbb F^n$ to
$\mathbb F^n,$ which are considered as right $\mathbb F$-vector spaces.
The corresponding unitary group over $\mathbb F$ is
\[ U(n,\mathbb F) = \{X\in M_n(\mathbb F): X^* X = I_n\},\]
where $X^* = \overline X^T$.
We denote by $\Delta$ the determinant on $M_q(\mathbb F),$ which is the
usual one for $\FF = \RR$ or $\CC$
and the
Dieudonn\'{e} determinant for $\FF = \HH$, i.e. $\Delta(X) = (\det_{\CC}
(X))^{1/2}$ when $X$ is regarded as a complex matrix of double size.

In this section, we consider the compact Grassmannians
$G_{p+q,q}(\mathbb F)= U/K$ with
$U= SU(p+q,\FF)$ and $K=S(U(p,\FF)\times U(q,\mathbb F)).$
Thus $U=SO(p+q)$, $SU(p+q)$ or $Sp(p+q)$ and $K=S(O(p) \times O(q))$,
$S(U(p) \times U(q))$ or $Sp(p) \times Sp(q)$, respectively.  We exclude the
case $p=q$ and assume that $p > q \ge 1$. 

Note that $SU(p+q)$ and $Sp(p+q)$ are simply connected, but $SO(p+q)$ is not,
nor is the Grassmannian $G_{p+q,q}(\mathbb R).$ So the general theory of Section
\ref{Symm_Raum} cannot be directly applied in the real case.
However,
$G_{p+q,q}(\mathbb R)$
has a simply connected double cover, namely 
 $SO(p+q)/SO(p)\times SO(q)$. This is just the Grassmannian of oriented
$p$-dimensional subspaces of $\mathbb R^{p+q}.$ Note also that $SO(p)\times
SO(q)$ is connected. Put
\[ K^\prime := \begin{cases} K & \text{if } \,\mathbb F = \mathbb C, \mathbb H\\
                SO(p)\times SO(q)&\text{if } \,\mathbb F = \mathbb R.
               \end{cases}\]
Then $K^\prime$ is connected and the Grassmannian $U/K^\prime$ is known to
be simply connected (c.f. \cite{Hel78}, Chap. X, Par.2).

We may therefore apply the theory  of Section \ref{Symm_Raum} to
$(U,K^\prime).$
For this, choose for the maximal abelian subspace $\lieb \subset \lieu$ the set
of all matrices $H_{x} \in M_{p+q} (\FF)$ of the form
\[ H_{x} = \begin{pmatrix} 0_{p\times p}& \begin{matrix} \underline x \\
                            0_{(p-q)\times q}\end{matrix}\\
\begin{matrix} \,  - \underline x& 0_{q\times(p-q)}\end{matrix}& 0_{q\times q}\end{pmatrix},
\]
where $\underline x := \textnormal{diag} (x_1, \ldots, x_q)$ is the $q \times q$ diagonal matrix corresponding to $x = (x_1, \ldots x_q) \in \RR^q$. See eg. \cite{Hel78}, p. 452 ff.

To keep notions convenient and consistent with \cite{Ro09}, we shall identify
the space $\liea= i\lieb$ (as
well as its dual $\liea^*$) with $\RR^q$, via $iH_x\mapsto x$.  In case $\mathbb
F=\mathbb R$, this amounts to an implicit use of the following isomorphism from
$\lieg = \liek \oplus \liep$ onto $\mathfrak{so}(p,q)$ which is described in
\cite{Hel78}, p. 453:

\[
\begin{pmatrix} A & iX \\ -iX^T & B \end{pmatrix} \mapsto \begin{pmatrix} A & X \\ X^T & B \end{pmatrix}
\]
The corresponding root system $\Sigma$ is of type $BC_q$ for $\mathbb F\in
\{\mathbb C, \mathbb H\}$ and of type $B_q$ for $\mathbb F = \mathbb R$. 
With our identification of $\liea$ and $\liea^*$ with $\mathbb R^q$, the
canonical choice of positive roots is $\Sigma^+ (BC_q) = \{ e_i, 2e_i \, : \, 1
\le i \le q \} \cup \{e_i \pm e_j \, : \,  1\le i < j \le q \}$ in the complex
and quaternionic case, and $\Sigma^+ (B_q) = \{ e_i \, : \, 1 \le i \le q \}
\cup \{e_i \pm e_j \, : \, 1\le i < j \le q \}$ in the real case. Here the $e_i$
are the standard basis vectors of $\RR ^q$. The Weyl group $W=W(\Sigma)$
is the hyperoctahedral group in all cases, and the Weyl chamber associated to
$\Sigma^+$ is given by
\[
\liea^+ = \{ x=(x_1, \ldots, x_q ) \in \RR^q: \, x_1 > x_2 > \ldots>  x_q > 0
\}.
\]
The roots with their multiplicities $m_{\alpha}$ are
given in the following table; the multiplicities depend on $p$, $q$ and the real
dimension $d=1,2,4$ of $\FF = \RR, \CC, \HH$. 

\begin{center}
\begin{tabular}{|c|c|}
\hline
root $\alpha$ & multiplicity $m_{\alpha}$ \\
\hline
$\pm e_i, \; \; 1 \le i \le q$ & $d(p-q)$ \\
\hline
$\pm 2e_i, \; \; 1 \le i \le q$ & $d-1$ \\
\hline
$\pm e_i \pm e_j, \; \; 1 \le i < j \le q$ & $d$ \\
\hline
\end{tabular}             
\end{center}

We will use the notation $m=(m_1,m_2,m_3)$ where $m_i$  denotes the
multiplicity on $\pm e_i$, $\pm 2 e_i$ or $\pm e_i \pm e_j$, respectively. The
Heckman-Opdam hypergeometric function of type $B_q$ coincides with a
hypergeometric function of type $BC_q$ having $m_2=0.$
For $BC_q$, the weight lattice is $\,  2\ZZ^q$ and the set of dominant weights
is given by 
\[
P^+ (BC_q)= \{ \lambda = (\lambda_1, \ldots, \lambda_q ) \in 2 (\ZZ^+)^q \, : \,
\lambda_1
\ge \lambda_2 \ge \ldots \ge \lambda_q \}.
\]
The closed fundamental alcove $A_0$ is given by
\[ A_0(BC_q) = \{ x\in \RR^q: \frac{\pi}{2} \geq x_1 \geq x_2 \geq \ldots \geq x_q\geq 0\}.\]
For root system $B_q$, the set of dominant weights is 
\[
P^+(B_q) =  \{ \lambda \in (\ZZ^+)^q \, : \, \lambda_1 \ge \lambda_2 \ge \ldots \ge \lambda_q; \; \text{ all } \lambda_i \text{ have same parity} \}
\]
and the alcove $A_0$ is bigger than in the $BC_q$-case.  It will however not be
needed in the sequel.

\medskip
The next theorem gives a  Cartan decomposition of $U$. It involves $K$ instead
of $K^\prime$ and the
$BC_q$-alcove  also in the real
case.

\begin{theorem}\label{KAK}
Let $U=SU(p+q, \FF)$ and $K=S(U(p,\FF)\times U(q,\mathbb F)).$  The group $U$
decomposes as $U=K S K$, where
\[
S = \left\{ b_x = \begin{pmatrix} \cos \underline x & 0_{q\times (p-q)} &  - \sin \underline x\\
0_{(p-q)\times q}& I_{p-q}& 0_{(p-q)\times q} \\
 \sin \underline x & 0_{q\times (p-q)} & \cos \underline x
            \end{pmatrix} \, : \, x \in A_0 \right\}
\]
with $A_0=A_0(BC_q).\,$ 
Every $u \in U$ can be written as $u=kb_xk'$ with $k,k' \in K$ and a unique $b_x \in S$. 
\end{theorem}

\begin{proof}
In the cases $\FF = \CC, \HH$ the group $U$ is simply connected and the result follows from Theorem 8.6 in Chapter VII of \cite{Hel78}: Put $\overline{Q_0} := \{H_{x} \, : \, x \in A_0 \}$. 
Then  a short calculation shows that $S = \exp \overline{Q_0}$.

In the case of $SO(p+q)$ this decomposition is explicitly given in
\cite{Vilenkin}, Section 15.1.9.
\end{proof}

As a consequence of this theorem, the  double coset space
$\,U//K = \{ KxK: x\in U\}\,$
is homeomorphic to the $BC_q$-alcove $A_0$ via $\,
Kb_xK\mapsto x\,.$

\medskip
Now we turn our attention to the spherical functions of $(U, K).$
Our first aim is to make the product formula
\begin{equation}\label{prodformula_abstract}
\psi(g) \psi(h) = \int_K \psi(gkh) dk
\end{equation}
explicit. For this, we may follow the argumentation of \cite{Ro09} (Section 2) in the non-compact dual cases. 
Since spherical functions on $U=KSK$ are $K$-biinvariant they are determined by their values on  $S$.
We consider 
\[ 
g:= \begin{pmatrix} u & 0 \\
0 & v \end{pmatrix} b_x \begin{pmatrix} \widetilde u & 0\\ 0 & \widetilde v
\end{pmatrix} \in KSK.
\]
and write $g$ in  $p\times q$ block notation as
\[ 
g = \begin{pmatrix} A(g) & B(g)\\
C(g) & D(g)
\end{pmatrix}.
\]
A short calculation then gives
\begin{equation}\label{D(g)}
D(g) = v \cos \xdiag \,\widetilde v
\end{equation}
where $\cos \xdiag \,= \text{diag}(\cos x_1, \ldots, \cos x_q)$. Note that $\cos x_i \in [0,1]$ since $x \in A_0$. 
We denote by $\text{spec}_s(X)$ the singular spectrum of $X \in M_q(\FF)$, that is
\[
\text{spec}_s(X) = \sqrt{\text{spec}(X^*X)} = (\sigma_1,\ldots, \sigma_q)\in \RR^q,
\]
with the singular values $\sigma_i$ of $X$  ordered by size:
$\sigma_1 \geq \ldots \geq \sigma_q\geq 0$.
Equation (\ref{D(g)}) implies that
the singular spectrum of $D(g)$ is given by \[\singspec (D(g)) = (\cos x_1,
\ldots, \cos x_q) =: \cos  x.\]
By our choice of the fundamental alcove $A_0$, we therefore have
\begin{equation}\label{x_von_D}
x = \arccos(\singspec(D(g)) \quad \forall g \in K b_x K, \, \,  x \in A_0,
\end{equation}
where $\arccos$ is also taken componentwise.
In order to evaluate formula \eqref{prodformula_abstract} explicitly, we write 
$b_x \in S$ in $p \times q$ block notation:
\[
b_x = \begin{pmatrix} A_x & B_x\\
C_x & D_x
\end{pmatrix}.
\]
Then for $\,k= \begin{pmatrix} u & 0\\
0 & v
\end{pmatrix} \in K\,$ we obtain by a short calculation  that
\[ 
D(b_xkb_y) = \, - \sin \,\underline x \,\sigma_0^* u \sigma_0 \sin \underline y \,+\,\cos \underline x \,v \cos \underline y
\] 
with the $p\times q$ block matrix
\[ 
\sigma_0 := \begin{pmatrix} I_q\\ 0\end{pmatrix}.
\]

Now let $\psi$ be a spherical function on $U$ and put $\, \widetilde\psi(x) := \psi(b_x)\,$ for $x \in A_0$. From  (\ref{x_von_D}) it follows  that $\widetilde\psi$ satisfies
\begin{equation}
\widetilde \psi(x) \widetilde\psi(y) \, = \, \int_K \widetilde\psi \left(\arccos \left( \text{spec}_s (D(b_x kb_y) \right) \right) dk.
\end{equation}

For our later extension of this product formula beyond the geometric cases, it
is important to rewrite it in a way where
the parameter $p$ is no longer contained in the domain of integration. Under the
technical assumption $p\geq 2q$, this can be done in  the same way as in
\cite{Ro09}, which leads to  the following

\begin{theorem}\label{prod_spherical}
Suppose that $p \ge 2q$. Define
\[D_q := \{ w \in M_q(\FF) : w^* w < I \},\quad 
 \gamma := d(q-\frac{1}{2})+1,\]
and for  $\mu \in \mathbb C$ with $\Re \, \mu > \gamma - 1$, put
\[
\kappa_{\mu} := \int_{D_q} \Delta (I - w^*w)^{\mu - \gamma} dw
\]
where $\Delta$ is the  determinant on $M_q(\FF).$ 
Then the spherical functions $\widetilde \psi (x) = \psi (b_x)$ 
of the Grassmannian $U/K =  G_{p+q,q}(\mathbb F)$ satisfy the
product formula
\begin{align*}
\widetilde \psi(x) \widetilde \psi(y) \, = \, \frac{1}{\kappa_{pd/2}} \int_{D_q}
\int_{U_0(q,\FF)} & \widetilde \psi \left(\arccos \left( \textnormal{spec}_s (
- \sin \,\underline x \, w \sin \underline y \,+\,\cos \underline x \,v \cos
\underline y) \right) \right) \\
& \cdot  \Delta(I - w^{\ast} w)^{pd/2 - \gamma} dv dw \quad (\forall \,x,y \in
A_0).
\end{align*}
Here $U_0(q,\FF)$ stands for the
connected component of $U(q, \mathbb F).$
\end{theorem}

\medskip

We are now going to identify the spherical functions of $
G_{p+q,q}(\mathbb F) =U/K$  as Heckman-Opdam  polynomials of type $BC_q.$

First, we determine
the heighest weight spaces $\Lambda_K(U)$ for the spherical
representations of $U$, c.f. Section \ref{Symm_Raum}.
The case $\mathbb F=\mathbb R,$ where neither $U$ nor $U/K$
is simply connected, needs special care.

\noindent
\textbf{Case 1:} $\mathbb F= \mathbb C$ or $\mathbb H$. In this case  $U$ is
simply connected, and by the Cartan-Helgason theorem we have
\[\Lambda^+_K(U) = P^+(BC_q) = \{ \lambda \in 2 (\ZZ^+)^q \, : \, \lambda_1 \ge
\lambda_2 \ge \ldots \ge \lambda_q \}.\]

\noindent
\textbf{Case 2:} $\mathbb F= \mathbb R.$ Besides $U/K= SO(p+q) / S(O(p)
\times O(q))$ consider its
simply connected double cover $U/K^\prime = SO(p+q)/SO(p)\times SO(q).$
For the latter, we obtain 
\[\Lambda_{K^\prime}(U)= P^+(B_q) = \{ \lambda \in (\ZZ^+)^q:  \lambda_1
\ge \lambda_2 \ge \ldots \ge \lambda_q; \,
\text{ all } \lambda_i \text{ have same parity} \} 
\] 

 In
Section 6 of \cite{Str86}, the spherical representations for each highest weight
are constructed explicitly. Here also the Grassmannians $SO(p+q)/S(O(p)
\times O(q))$   are
considered. The fact that $K$ is larger than $K^\prime$ implies a further
invariance: all $\lambda_i$ have to be even and therefore
\[
\Lambda^+_K(U) = \{ \lambda \in 2 (\ZZ^+)^q \, : \, \lambda_1 \ge \lambda_2 \ge
\ldots \ge \lambda_q \}.
\]

So over all fields ($\RR$, $\CC$ and $\HH$), the spherical functions of
$(U,K)$ are indexed by the dominant weights of the root system $BC_q$. We
will use the notation
\begin{equation}
P^+ := P^+(BC_q) = \{ \lambda \in 2 (\ZZ^+)^q \, : \, \lambda_1 \ge \lambda_2
\ge \ldots \ge \lambda_q \}
\end{equation}
as well as $A_0 = A_0(BC_q)$ in the following.
From Theorem \ref{ident_C_H} and the above considerations, we conclude

\begin{theorem}\label{geometric}
 The spherical functions of the compact Grassmanian $ G_{p+q,q}(\mathbb
F) = U/K$, with $\FF \in \{\RR, \CC, \HH\},$ are indexed by $P^+$ and given by
Heckman-Opdam polynomials of type $BC_q$,
\begin{equation}\label{indent}
\widetilde \psi_{\lambda} (x) = \psi_{\lambda}(b_x) = F_{BC_q}(\lambda + \rho,m;
ix) =  R_\lambda(x) \quad (\lambda \in P^+),
\end{equation}
with the (geometric) multiplicity  $m= (d(p-q), d-1, d).$
\end{theorem}

Under the assumption $p \ge 2q$  Proposition \ref{prod_spherical} implies that
\begin{align*}
R_{\lambda} (x) R_{\lambda} (y) =   \frac{1}{\kappa_{pd/2}} 
& \cdot \int_{B_q} \int_{U_0(q,\FF)} R_{\lambda} \left(d(\underline x,
\underline y, v,w) \right) \Delta(I - w^{\ast} w)^{pd/2 - \gamma} dv dw
\end{align*}
for all $x,y \in A_0$, where
\begin{align}\label{d_formula}
d(\underline x, \underline y, v,w) := \arccos \left( \textnormal{spec}_s(- \sin \underline x \, w \sin \underline y + \cos \underline x \,v \cos \underline y)   \right).
\end{align}

The next step is analytic continuation. Fix $q$ and $d = \text{dim}_{\RR} \FF$.
For $\mu \in \CC$ with $\text{Re}\,\mu > \gamma -1$
and  $\lambda \in P^+$ consider the $BC_q$-Jacobi polynomials
\[
R_{\lambda}^{\mu}(x) := F_{BC_q} (\lambda + \rho_\mu, m_{\mu}; ix),
\]
with the multiplicity 
\begin{equation}\label{mmu}
m_{\mu} = (2 \mu - dq,d-1,d).
\end{equation}
Note that $\mu \mapsto
R_\lambda^{\mu}(x)$ is (by analytic extension) holomorphic on  $\{\text{Re}
\,\mu >\gamma -1\}.$
For $\mu = pd/2$ with $p \in \NN$ this gives the geometric cases of Theorem
\ref{geometric}.

\begin{theorem}\label{prod_Jacobi}
For $\mu \in \CC$ with $\Re \mu > \gamma - 1$ and $\lambda \in \Lambda^+$ the
Jacobi polynomials $R_\lambda^\mu$ satisfy the
product formula
\begin{align*}
R^{\mu}_{\lambda}(x) R^{\mu}_{\lambda}(y) =  \frac{1}{ \kappa_{\mu}} \int_{D_q}
\int_{U_0(q,\FF)} R^{\mu}_{\lambda} \left(d(\underline x, \underline y, v,w)
\right) \Delta(I - w^{\ast} w)^{\mu - \gamma} dv dw.
\end{align*}
\end{theorem}

\begin{proof}
The proof is a copy of the first part of the proof of Theorem 4.1 in
\cite{Ro09}.
Replace the $(R,k)$-notation by our $(\Sigma,m)$-notation and the product
formula by our product formula. Then rewrite the claimed formula in terms of
$P^{\mu}_{\lambda} := c(\lambda+\rho, m_{\mu})^{-1} \cdot R^{\mu}_{\lambda}$
(the standard Heckman-Opdam normalization):
\[
P^{\mu}_{\lambda}(x) P^{\mu}_{\lambda}(y) =
\frac{c(\lambda+\rho, m_{\mu})^{-1}}{ \kappa_{\mu}} \int_{D_q}
\int_{U_0(q,\FF)} P^{\mu}_{\lambda} \left(d(\underline x, \underline y, v,w)
\right) \Delta(I - w^{\ast} w)^{\mu - \gamma} dv dw
\]
For fixed $\lambda \in P^+$, the function $c(\lambda + \rho, m_{\mu})$ is
bounded away from zero as
$\mu \to \infty$ in the half plane $H=\{\mu \in \CC \, : \, \text{Re } \mu >
\gamma - 1 \}$ (see \cite{Ro09}). Then one uses the fact that the coefficients
of the $P_{\lambda}^\mu$ with respect to the exponential basis $\{ e^{i \nu}:
\, \nu \in P \}$ are rational, and that the integral
\[
\frac{1}{\kappa_{\mu}} \int_{D_q} | \Delta(I-w^* w)^{\mu - \gamma} | dw
\]
converges for $\text{Re } \mu > \gamma - 1$ and is of polynomial growth
as
$\mu \to \infty$ in $H$. This allows to apply Carlson's theorem. For details  we
refer to \cite{Ro09}.
\end{proof}

%
\section{The rank one case}\label{rg1}
%

At this point it is worthwile to see how our product formula
for Heckman-Opdam Jacobi polynomials generalizes the product formula of
classical one-variable Jacobi polynomials for certain indices.

The classical Jacobi polynomials with indices $\alpha, \beta >-1$ are given by
\begin{equation}\label{klass_Jacobi}
 P^{(\alpha, \beta)}_n (x) = \frac{(\alpha + 1)_n}{n!} \,{}_2 F_1 \bigl(\alpha + \beta+n+1, -n , \alpha + 1; \frac{1-x}{2} \bigr)
\end{equation}
where ${}_2F_1$ is the Gaussian hypergeometric function and
$(a)_n= \Gamma(a + n)/\Gamma(a).$ We renormalize:
\[ 
R_n^{(\alpha,\beta)}(x) := \frac{n!}{(\alpha + 1)_n}P_n^{(\alpha, \beta)}(x).
\]

Let us consider the Heckman-Opdam theory in the rank one case. The
root system is
$BC_1 = \{ \pm e_1, \pm 2 e_1 \}$ in $\liea \cong \RR$ and we
denote the multiplicity by $m=(m_1, m_2).$
According to the example in \cite{Op95}, p. 89f, the hypergeometric function
$F_{BC_1}$ is given by
\begin{equation}\label{hyper_ident}
F_{BC_1} (\lambda,m; x) = {}_2 F_1 \bigl(a,b,c; \frac{1}{2} ( 1 - \cosh
2x)\bigr)
\end{equation}
with
\begin{equation}\label{abc_Rang1}
a = \frac{1}{2} \bigl(\lambda + \frac{1}{2} m_1 + m_2 \bigr),
\quad b= \frac{1}{2} \bigl( - \lambda + \frac{1}{2} m_1 + m_2 \bigr)
\textnormal{ and } c=  \frac{1}{2} \bigl(1 + m_1 + m_2 \bigr).
\end{equation}

The dominant weights and the fundamental alcove  are
\[ P^+ =  2 \ZZ^+ , \quad A_0 = \big[0, \frac{\pi}{2}\big].\]
According to Lemma \ref{Hyper_Jacobi} we have
$\, F_{\lambda + \rho}(ix) =  R_{\lambda}(x)
\,$
where $\lambda = 2n \in 2\mathbb Z^+$ and $\rho = \frac{1}{2}m_1 + m_2\,.$
Equation (\ref{hyper_ident}) becomes
\begin{align*}
F_{BC_1} (\lambda,m; ix) &= {}_2 F_1 \bigl( n +\frac{1}{2} m_1 +  m_2, -n,
\frac{1}{2} (1 + m_1 + m_2);
\frac{1}{2} ( 1 - \cos 2x)\bigr) .
\end{align*}
In view of (\ref{klass_Jacobi}), we conclude
\begin{equation}\label{ident_polynome}
R_{\lambda} (x) = R_n^{(\alpha, \beta)} (\cos 2x)
\end{equation}
with
\begin{equation*}\label{alpha_beta}
\alpha = \frac{1}{2} (m_1 + m_2 -1), \quad \beta =
\frac{1}{2} (m_2 - 1).
\end{equation*}
In particular, we obtain the well-known fact that the spherical functions of the
rank one symmetric space $U/K=  G_{p+1,1}(\mathbb F)$ -- which is
just the $p$-dimensional projective space $\mathbb P^p(\mathbb F)$ --
are given in terms of classical Jacobi polynomials $R_n^{(\alpha, \beta)}$
with $\alpha = (dp-2)/2, \, \beta = (d-2)/2$.

In the real case $\FF = \RR,$ the
spherical functions of $\mathbb P^p(\mathbb R)$ are Gegenbauer
polynomials of even degree. In fact, \eqref{ident_polynome} becomes
\[
R_{\lambda} (x) = R_n^{(\alpha, - \frac 1 2)} (\cos 2x) = R_{2n}^{(\alpha,
\alpha)} (\cos x) \quad \text{ with }\,\alpha = \frac{1}{2} (p-2).\]
In the 1970ies, Koornwinder devoted a series of papers to the product formula
for one-variable  Jacobi polynomials, see e.g. \cite{Ko74}. For arbitrary
$\alpha > \beta > - \frac 1 2$, it is given by
\begin{align*}
R_n^{(\alpha, \beta)} (t) R_n^{(\alpha, \beta)} (s)= \int_0^1 \int_0^{\pi} & R_n^{(\alpha, \beta)}\bigl( \frac 1 2 (1+t)(1+s) + \frac 1 2 (1-t) (1-s) r^2 \\
& + \sqrt{1 - t^2} \sqrt{1-s^2}\, r \cos \theta -1 \bigr) \, dm_{\alpha,\beta}(r, \theta)
\end{align*}
with
\[
dm_{\alpha,\beta}(r, \theta) = c_{\alpha, \beta} (1-r^2)^{\alpha - \beta - 1} (r \sin \theta)^{2 \beta} r \, dr d\theta
\]
and
\[
\frac{1}{c_{\alpha,\beta}} = \int_0^1 \int_0^{\pi} (1-r^2)^{\alpha - \beta - 1} (r \sin \theta)^{2 \beta} r \, dr d\theta.
\]

Now consider the product formula from Theorem \ref{prod_Jacobi} for rank
$q=1$. Here $\gamma =\frac{d}{2} +1,$ and
we restrict to real parameters $\mu > \frac{d}{2}.$  Recall
again the identification \eqref{ident_polynome}. As $m_1 = 2\mu -d $ and $m_2 =
d-1,$ we have $\alpha = \mu-1$ and $\beta = (d-2)/2.$
The domains of integration reduce to $D_1 = \{ w \in \FF \, : \, |w| < 1\}$ and
$U_0(1) = \{ v \in \FF \, : \, |v| = 1 \}_0$. Furthermore,
\[
d(x,y,v,w) = \arccos | - w \sin x \sin y + v \cos x \cos y |.
\]
The $U_0(1)$-integral cancels under the coordinate transform $w' := v^{-1} w$.
Using $\cos 2x = 2 \cos^2 x - 1$ we obtain for $\alpha > \beta =\frac{d}{2}-1$
the product formula
\begin{align}\label{rg1_vorstufe}
 &R^{(\alpha, \beta)}_n (\cos 2x)  R^{(\alpha, \beta)}_n (\cos 2y)\,=\notag   \\
 &\,\,\,\,\frac{1}{\kappa_{\alpha +1}}
\int_{D_1}   R^{(\alpha, \beta)}_n (2 \, |- z \sin x \sin y + \cos x \cos y|^2
-1 )\cdot (1-|z|^2)^{\alpha -d/2} dz.
\end{align}
Let us sketch the further calculations only in the case $\FF = \CC$, where
$d=2.$ We introduce polar coordinates $z=re^{i \theta}$ and put
$t := \cos 2x$, $s:= \cos 2y$. Then use the identities $\sin^2 x = \frac 1 2
(1-t)$, $\sin x \cos x = \frac 1 2 \sqrt{1-t^2}$ and $\cos^2 x = \frac 1 2
(1+t)$. The constant $\kappa_{\alpha +1}$ is given by
\[
\kappa_{\alpha +1} = 2 \pi \int_0^1 (1-r^2)^{\alpha -1} r \, dr =
\frac{\pi}{\alpha}.
\]
We conclude from (\ref{rg1_vorstufe}) exactly the product formula for the Jacobi
polynomials $R_n^{(\alpha, \beta)}$ with $\alpha >0, \,\beta = 0$:
\begin{align*}
R_n^{(\alpha, 0)} (t) R_n^{(\alpha, 0)} (s) = & \frac{2\alpha}{\pi} \int_0^1
\int_0^{\pi}  R_n^{(\alpha, 0)} \Bigl( \frac 1 2 (1+t)(1+s) + \frac 1 2 (1-t)
(1-s) r^2 \notag \\
& + \sqrt{1 - t^2} \sqrt{1-s^2} \, \, r \cos \theta -1 \Bigr) (1-r^2)^{\alpha
-1} r \, dr d\theta.
\end{align*}

%
\section{Hypergroup structures on the alcove}\label{sec_hyper}
%

In this section we shall see that the product formula of Theorem
\ref{prod_Jacobi} leads to three
continuous series (for $d=1,2,4$) of positivity-preserving convolution algebras
on the fundamental alcove $A_0 = \{ x \in \RR^q \, : \, \frac{\pi}{2} \ge
x_1 \ge \ldots \ge x_q \ge 0 \},$ which are compact commutative
hypergroups with normalized Jacobi polynomials as characters. In the geometric
cases ($\mu = pd/2$), these hypergroup
convolutions are just given by the double coset convolutions on a double coset
space $U//K$ which may be identified with $A_0$ according to Theorem \ref{KAK}.

To start with, let us briefly recall some basics from hypergroup theory. For a detailed treatment, the reader is referred to \cite{Jew75}. Hypergroups generalize the convolution algebras of locally compact groups, with the convolution product of two point measures $\delta_x$ and $\delta_y$ being in general not a point measure again but a probability measure with compact support depending on $x$ and $y$.
\begin{definition}
A hypergroup is a locally compact Hausdorff space $X$ with a weakly continuous, associative convolution $\ast$ on the space $M_b(X)$ of regular bounded Borel measures on $X$, satisfying the following properties:
\begin{enumerate}
\item The convolution product $\delta_x \ast \delta_y$ of two point measures is a compactly supported probability measure on $X$, and $\text{supp}(\delta_x \ast \delta_y)$ depends continuously on $x$ and $y$ with respect to the so-called Michael topology on the space of compact subsets of $X$ (see \cite{Jew75}).
\item There exists a (necessarily unique) neutral element $e \in X$ satisfying $\delta_e \ast \delta_x = \delta_x \ast \delta_e = \delta_x$ for all $x \in X$.
\item There exists a (necessarily unique) continuous involution $x \mapsto
\overline x$ on $X$  such
that $\delta_{\overline x} \ast \delta_{\overline y} = (\delta_y \ast
\delta_x)^-$ and $x = \overline y \iff e \in \textnormal{supp}(\delta_x \ast
\delta_y)$. Here the measure $\mu^-$ is given by $\mu^- (A) = \mu(\overline
A)$.
\end{enumerate}
The hypergroup is called commutative if the convolution is commutative.
\end{definition}

Note that due to weak continuity, the convolution of measures on a hypergroup is
uniquely determined by the convolution of point measures.

Every commutative hypergroup $X$ has a unique (up to a multiplicative factor)
Haar measure $\omega$,
that is a positive Radon measure with the property
\[
\int_X f(x \ast y) d\omega(y) = \int_X f(y) d\omega(y) \quad (\forall x \in X, f \in C_c(X)),
\]
where we use the notation $f(x \ast y) := (\delta_x \ast \delta_y)(f)$.

The dual space of a commutative hypergroup $X$ is defined by
\[
\widehat X := \{ \phi \in C_b(X) \, : \; \phi \ne 0, \; \phi(\overline x) = \overline{ \phi (x)}  \, \textnormal{ and } \, \phi (x \ast y) = \phi(x) \phi(y) \}.
\]
The elements of $\widehat X$ are called characters of $X$. The dual of a
commutative hypergroup is a locally compact Hausdorff space with the topology of
locally uniform convergence. In the case of a compact hypergroup $X$ the dual
$\widehat X$ is discrete. The Fourier transform on $L^1(X, \omega)$ is defined
by
\[\widehat f (\phi) := \int_X f(x) \overline{\phi(x)} d\omega(x), \,\, \phi \in
\widehat X.\]
 It is injective and there exists a unique positive Radon measure
$\pi$ on $\widehat X$, called the Plancherel measure of $X$, such that $f
\mapsto \widehat f$ establishes an isometric isomorphism from $L^2(X, \omega)$
onto $L^2(\widehat X, \pi)$.

\begin{example} (Double coset hypergroups) \\
Let $G$ be a locally compact group with compact subgroup $K$ and denote by $dk$
the normalized Haar measure on $K$. Then there is a natural hypergroup structure
on the set of double cosets $G//K = \{ KxK \, : \, x \in G \}$ which is given by
\[
\delta_{KxK} \ast \delta_{KyK} = \int_{K} \delta_{KxkyK} \, dk, \quad x,y \in G.
\]
The neutral element is $K=KeK$ and the involution is given by $(KxK)^- = Kx^{-1}K$ (see Theorem 8.2B in \cite{Jew75}). The double coset hypergroup $(G//K, \ast)$ is commutative if and only if $(G,K)$ is a Gelfand pair. 
\end{example}

We now return to the setting of Theorem \ref{prod_Jacobi}. 

\begin{theorem}\label{convo}
Let $\mu \in \mathbb R$ with $\mu > \gamma-1$. Then the
probability measures
\begin{align*}
\delta_{x} \ast_{\mu} \delta_{y} (f) := \frac{1}{\kappa_{\mu}} \int_{D_q}
\int_{U_0(q,\FF)} f \left( d(\underline x, \underline y, v,w) \right) \Delta(I -
w^{\ast} w)^{\mu - \gamma} dv dw
\end{align*}
with
\[
d(\underline x, \underline y, v,w) = \arccos \left( \textnormal{spec}_s(- \sin \underline x \,w \sin \underline y + \cos \underline x \, v \cos \underline y)    \right)
\]
for $x,y \in A_0$ define a commutative hypergroup structure on the compact alcove $A_0$. The neutral element is $0$ and the involution is the identity mapping.
\end{theorem}

Note that in the geometric cases $\mu =pd/2,$ the convolution $*_\mu$ on $A_0$
is just the convolution of the corresponding double coset hypergroup $U//K$.

\begin{proof} We use standard arguments. First, 
the integral defining the convolution is invariant under
$v \mapsto v^*$, $w \mapsto w^*$ and $d(\underline x, \underline y, v,w) =
d(\underline y, \underline x, v^*,w^*)$. Therefore $\ast_\mu$ is
commutative. For associativity let $x,y,z \in A_0$. Then for $f\in C(A_0),$
\begin{align*}
\delta_{x} \ast_{\mu} ( \delta_{y} \ast_{\mu} \delta_{z}) (f) =&
\frac{1}{\kappa_{\mu}^2}
\int_{B_q \times U_0(q)} \int_{B_q \times U_0(q)}
f(D(x,y,z,v,w,v',w'))\cdot & \\ & \cdot\Delta(I-w^*w)^{\mu-\gamma}
\Delta(I-(w')^*w')^{\mu-\gamma} dvdwdv'dw' & =: I(\mu)
\end{align*}
with a certain $A_0$-valued argument $D$, which is independent of $\mu$.
The same is true for
\[
( \delta_{x} \ast_{\mu}  \delta_{y}) \ast_{\mu} \delta_{z} (f) =: I'(\mu)
\]
with a $\mu$-independent argument $D'$ instead of $D$. The integrals $I(\mu)$
and $I'(\mu)$ are well
defined and holomorphic in $\{ \mu \in \CC \, : \, \Re \mu > \gamma -1 \}$.
The convolution is associative in the geometric cases $\mu = pd/2$.
Analytic continuation then yields associativity for all $\mu$ with $\Re \mu >
\gamma
-1$ as in \cite{Ro07}. 
Weak continuity of the convolution  follows from the continuity of the
mapping $(x,y,v,w) \mapsto f(d(x,y,v,w))$ on $A_0^2 \times B_q \times U_0(q)$.
It is also
obvious that $0$ is neutral. So only the support continuity 
and the fact that the identity mapping is a hypergroup
involution remain. As the support of $\delta_{x} \ast_{\mu} \delta_{y} $ is
independent of $\mu$, 
it suffices to verify both statements  in the geometric cases $U//K$. But these
are
known to correspond to double coset hypergroups, which immediately implies the
support continuity. In the geometric cases,  the involution is induced by
the group inversion on $U$, and hence by the mapping $x \mapsto -x$ on
$\RR^q\cong \liea.$ A short calculation shows that $b_{-x} \in Kb_{x}K$ on
$U//K$ and
therefore the involution on $U//K$ is the identity.
In fact,
\begin{align*}
& \begin{pmatrix} \cos \underline x & &  \sin \underline x \,\\ & I_{p-q} &
\\ -\sin \underline x & & \cos \underline x \,\end{pmatrix} =\\
&\qquad =\begin{pmatrix} - I_q & & \\ & I_{p-2q,q} & \\ & & I_q \end{pmatrix}
\begin{pmatrix} \cos \underline x & & - \sin \underline x \\ & I_{p-q} &
\\ \sin \underline x & & \cos \underline x \end{pmatrix} \begin{pmatrix}
- I_q & & \\ & I_{p-2q,q} & \\ & & I_q \end{pmatrix}
\end{align*}
where $I_{n,m} = \text{diag}(1,\ldots,1,-1, \ldots,-1)$ denotes the
diagonal matrix with the first $n$ entries equal to $1$ and the last $m$
entries equal to $-1$. Then \\$\textnormal{det} \begin{pmatrix} - I_q & \\
& I_{p-2q,q} \end{pmatrix} = 1$.
\end{proof}

\begin{proposition}
The support of $\delta_x \ast_{\mu} \delta_y$ satisfies
\[
\textnormal{supp}(\delta_x \ast_{\mu} \delta_y) \subseteq \{ z \in A_0 \, :
\, \| z \|_{\infty} \le \| x \|_{\infty} + \| y \|_{\infty} \}
\]
where $\| \cdot \|_{\infty}$ is the maximum norm in $\RR^q$.
\end{proposition}

\begin{proof} This is more involved than the corresponding statement in the non-compact case (\cite{Ro09}).
For a matrix $A\in M_q(\mathbb F)$ we denote again by
\[\textnormal{spec}_s (A)
= (\sigma_1(A), \ldots, \sigma_q(A))\in \mathbb R^q\]
 the singular values of
$A$, decreasingly ordered by size.
Write $\|A\| =  \| \singspec (A) \|_{\infty} = \sigma_1 (A)$ for the spectral
norm of $A$. We need  the following estimates from Theorem 3.3.16 in
\cite{HJ91}:
\begin{eqnarray}
|\sigma_q (A+B) -  \sigma_q (A) | \le \sigma_1 (B) \label{est_1} \\
\sigma_q(AB) \le \sigma_q (A) \sigma_1(B) \label{est_2}
\end{eqnarray}
These estimates are only stated for $\FF= \RR, \CC$ but in the case of a
quaternionic matrix $A \in M_q (\HH)$  we simply consider the corresponding
complex matrix $\chi_A \in M_{2q} (\CC)$, namely
\[
\chi_A = \begin{pmatrix} A_1 & A_2 \\ - \overline{A_2} & \overline{A_1} \end{pmatrix}
\]
where $A_1, A_2$ are complex $q \times q$-matrices such that $A = A_1 + A_2 j$. The map $M_q (\HH) \to M_{2q} (\CC)$, $A \mapsto \chi_A$ is a homomorphism and $\chi_{A^*} = \left( \chi_A \right)^*$. Moreover, $\textnormal{spec}_s (A) = \textnormal{spec}_s (\chi_A)$ where in the second set each singular value appears twice (see \cite{Zhang} for a survey about quaternionic matrices).

Let $\xi := \cos \underline x \, v \, \cos \underline y - \sin \underline x \, w \, \sin \underline y$. By \eqref{est_1},
\[
\sigma_q (\xi) \ge \sigma_q (\cos \underline x \, v \, \cos \underline y) - \sigma_1 (\sin \underline x \, w \, \sin \underline y).
\]
Since $\sin x$ is increasing on $[0, \pi/2]$ we get (using submultiplicativity)
\[
\sigma_1 (\sin \underline x \, w \, \sin \underline y) = \| \sin \underline x \, w \, \sin \underline y \| \le \| \sin \underline x \| \| \sin \underline y \| = \sin \| x \|_{\infty} \sin \| y \|_{\infty}.
\]
On the other hand, if  $\cos y_i \ne 0$ for all $i$, then by \eqref{est_2}  
\[
\sigma_q (\cos \underline x \, v \, \cos \underline y) \ge\, \frac{\sigma_q (\cos \underline x \, v ) }{\sigma_1 \left( (\cos \underline y)^{-1} \right)} \ge \cos \| x \|_{\infty} \cos \| y\|_{\infty}.
\]
Therefore
\[
\sigma_q (\xi) \ge \cos \| x \|_{\infty} \cos \| y\|_{\infty} - \sin \| x \|_{\infty} \sin \| y \|_{\infty} = \cos (\|x\|_{\infty} + \| y \|_{\infty})
\]
This implies the claim, because $\arccos $ is decreasing. If $\cos y_i = 0$ for some $i$, the estimate remains by continuity since the eigenvalues of a matrix depend continuously upon its entries; see e.g. \cite{HJ91}, p. 396.
\end{proof}

Because of Theorem \ref{prod_Jacobi} the Jacobi polynomials
$R_\lambda := R_{\lambda}^{\mu}$ are multiplicative, 
\begin{equation}\label{R_mult}
R_{\lambda}(x) R_{\lambda}(y) = R_{\lambda}(x \ast_{\mu} y).
\end{equation}

\begin{lemma}\label{La_R_reell}
Assume that the Weyl group $W= W(\Sigma)$ contains the reflection $\sigma : x
\mapsto -x$. Then for nonnegative multiplicities, the associated Heckman-Opdam
polynomials $P_{\lambda}$ are real-valued on $\RR^q$. In particular, this holds
for the root systems $\Sigma = B_q$, $C_q$ and $BC_q$.
\end{lemma}

\begin{proof}
This is immediate from identity \eqref{conjugation}.\end{proof}

In our situation, the Jacobi polynomials
$R_{\lambda} = R_\lambda^\mu,  \, \lambda \in P^+$ are therefore
indeed characters of the hypergroup $(A_0,\ast_\mu)$. It is part of the
following theorem that they make up the complete dual.

\begin{theorem}
\begin{enumerate}
 \item[\rm{(a)}] The Haar measure of the hypergroup $(A_0, \ast_\mu)$ is 
\[d \omega(x) = w_m(x)dx = \prod_{\alpha \in \Sigma^+} \left| e^{i\langle
\alpha,x\rangle} -
e^{- i\langle \alpha, x \rangle} \right|^{ m_\alpha} \; dx\]
where $m = m_\mu$ as defined in \eqref{mmu}.
 \item[\rm{(b)}] The dual space is $(A_0, *_\mu)^\wedge\, = \{ R_{\lambda} \, :
\, \lambda \in P^+ \}$.
\end{enumerate}
\end{theorem}

\begin{proof} (a)\enskip
For  $R_{\lambda}$ with $\lambda \ne 0$ we have
$\, \int_{A_0} R_{\lambda} d\omega = 0
\,$ 
since $R_\lambda$ is orthogonal to $R_0 = 1.$
In view of (\ref{R_mult}), we obtain 
\begin{align*}
\int_{A_0} R_{\lambda}(x \ast_\mu y)  d\omega(y) = R_{\lambda} (x) \int_{A_0} R_{\lambda} (y) d\omega(y) = 0.
\end{align*}
By linearity, the above equation holds for all $W$-invariant trigonometric
polynomials. By the Stone-Weierstrass theorem, $\mathcal T^W$ is $\| \cdot
\|_{\infty}$-dense in $C(A_0).$ Now the
assertion follows from the $    \| \cdot \|_{\infty}$-continuity of the
hypergroup translation (see Lemma 3.3B in \cite{Jew75}).

(b)\enskip
We already know that the $R_{\lambda}$ are characters of our hypergroup. In
general, the characters of a compact commutative hypergroup $X$ form an
orthogonal basis of $L^2(X, d\omega)$. The proof is the same as in the case of a
compact group and uses the Plancherel Theorem, see Theorem 3.5 in \cite{Dun73}.
The Jacobi polynomials form already an orthogonal basis of $L^2(A_0, \omega).$
So there are no additional characters.
\end{proof}

\begin{remark}
For a general commutative hypergroup $X$ the set of bounded semi-characters 
\[
\chi_b (X) := \{ \phi \in C_b(X) \, : \; \phi \ne 0 \, \textnormal{ and } \, \phi (x \ast y) = \phi(x) \phi(y) \}
\]
may not coincide with the dual $\widehat X$. However, if $X$ is compact
(more general: of subexponential growth), then it can be shown by Banach
algebraic methods that
$\widehat X = \chi_b(X)$; see Theorem 2.5.12 in \cite{BH95}. But Lemma
\ref{La_R_reell}, which leads to this identity in our present case, is also of
some interest in its own.

\end{remark}

We identify the dual of the hypergroup $(A_0, \ast_\mu)$ with the set of
dominant weights via the mapping $(A_0)^{\wedge} \to P^+$, $R_{\lambda}
\mapsto \lambda$. 

\begin{proposition}\label{Plancherel}
The Plancherel measure of the hypergroup $(A_0, \ast_{\mu})$ is the following
measure on $P^+$:
\[
\pi = \sum_{\lambda \in P^+} r_{\lambda} \delta_{\lambda}
\]
with
\[r_{\lambda} := \bigl(\int_{A_0} |R_\lambda|^2 d\omega\bigr)^{-1}.\]
\end{proposition}

\begin{proof}
The set  $\{ \sqrt{r_{\lambda}} R_{\lambda} \, : \, \lambda \in P^+ \}$ is an
orthonormal basis of $L^2(A_0, \omega)$. Thus for $f\in L^2(A_0,\omega),$
\[
\int_{A_0} |f|^2 d\omega = \sum_{\lambda \in P^+} r_{\lambda}  |\langle f,
R_{\lambda} \rangle |^2 = \sum_{\lambda \in P^+}  r_{\lambda} |\widehat
f(\lambda)|^2 = \int_{P^+} |\widehat f|^2 d\pi.
\]
\end{proof}

\end{document}